\documentclass[11pt]{amsart}
\usepackage{geometry}                
\usepackage[parfill]{parskip}    
\usepackage{graphicx}
\usepackage{amssymb}
\usepackage{subfigure}
\usepackage{amsmath}
\usepackage{color}
\usepackage{setspace}

\usepackage{pdfsync}
\swapnumbers 
\numberwithin{equation}{section} 
\theoremstyle{plain} 
\newtheorem{thm}{Theorem}[section] 
\newtheorem{cor}[thm]{Corollary} 
 
\newtheorem{lem}[thm]{Lemma} 
 
\newtheorem{prop}[thm]{Proposition}

\theoremstyle{definition} 
\newtheorem{defin}[thm]{Definition}

\newtheorem{mex}[thm]{Main Example}

\def\pp{\medskip{\parindent 0pt \it Proof.\ }} 

\def\Z{{\mathbb Z}}

\def\arr{\protect\operatorname{arr}}
\begin{document}
\title[A Generalization of Turaev's Virtual String Cobracket]{A Generalization of Turaev's Virtual String Cobracket and Self-Intersections of Virtual Strings}
\author{Patricia Cahn}\address{Patricia Cahn, Department of Mathematics, University of Pennsylvania, David Rittenhouse Lab. 209 South 33rd Street, 
Philadelphia, PA 19104-6395, USA} \email{pcahn@sas.upenn.edu}

\begin{abstract}  Previously we defined an operation $\mu$ that generalizes Turaev's cobracket for loops on a surface.  We showed that, in contrast to the cobracket, this operation gives a formula for the minimum number of self-intersections of a loop in a given free homotopy class.  In this paper we consider the corresponding question for virtual strings, and conjecture that $\mu$ gives a formula for the minimum number of self-intersection points of a virtual string in a given virtual homotopy class. To support the conjecture, we show that $\mu$ gives a bound on the minimal self-intersection number of a virtual string which is stronger than a bound given by Turaev's virtual string cobracket.  We also use Turaev's based matrices to describe a large set of strings $\alpha$ such that $\mu$ gives a formula for the minimal self-intersection number $\alpha$.  Finally, we construct an example that shows the bound on the minimal self-intersection number given by $\mu$ is always at least as good as, and sometimes stronger than, the bound $\rho$ given by Turaev's based matrix invariant. 
\end{abstract}
\maketitle
\section{Introduction}
The goal of this paper is to estimate, and in many cases compute precisely, the minimum number of double points of a flat virtual knot in a given virtual homotopy class. To do this we define a generalization of Turaev's Lie cobracket on the vector space generated by homotopy classes of flat virtual knots. We then compare our estimates to estimates given by Turaev's Lie cobracket, and Turaev's based matrix invariant.  We conjecture that our generalization of Turaev's cobracket always computes the minimum number of double points.

 {\bf Notation and conventions.} A virtual string, or flat virtual knot, is a combinatorial generalization of a curve on a surface.  Virtual knots were introduced by Kauffman \cite{Kauffman}.  In this paper we use Turaev's terminology, and refer to flat virtual knots as virtual strings \cite{Turaev2}.  We represent virtual strings with Gauss diagrams, following the conventions in \cite{Turaev2}. The {\it minimal self-intersection number} of a virtual string $\alpha$ is the minimum number of arrows of a string in the virtual homotopy class of $\alpha$.  We refer the reader to Section \ref{sec:strings} for precise definitions of virtual strings and their homotopy classes.  For a (reduced) finite linear combination $L=\sum_i a_i e_i$ in the free $\mathbb{Z}$-module $\Z[S]$, where $a_i\in\Z$ and $e_i\in S$, we put $t(L)=\sum_i|a_i|$ and call it the {\it number of terms of $L$}.\\

Turaev \cite{Turaev2} defined a Lie cobracket $\nu$ on the free $\Z$-module generated by the set of nontrivial homotopy classes of virtual strings.  We define an operation $\mu$ on the same $\Z$-module.  Turaev's cobracket factors through $\mu$, so we view $\mu$ as a generalization of $\nu$.  Both operations give lower bounds on the minimal self-intersection number of a virtual string $\alpha$.  Let $[\alpha]$ denote the homotopy class of the virtual string $\alpha$.  Let $m([\alpha])$ denote the minimal self-intersection number of $\alpha$.  The bounds on $m([\alpha])$ given by $\mu$ and $\nu$ are
\begin{equation}\label{virtualmubound}
m([\alpha])\geq t(\mu([\alpha]))/2 +n-1, 
\end{equation}
and
\begin{equation}\label{virtualnubound}
m([\alpha])\geq t(\nu([\alpha]))/2 +n-1,
\end{equation}
where $n\geq 1$ is the largest integer such that a representative $\beta \in [\alpha]$ with minimal self-intersection in $[\alpha]$ can be realized as a curve $B$ on a surface $F$, and $\langle B\rangle=\langle \gamma \rangle ^n$ for some $\langle \gamma \rangle$ in $\pi_1(F)$.   (Recall that a virtual string $\alpha$ is {\it realized} as a curve $A$ on an oriented surface $F$ if the Gauss diagram of $A$ is $\alpha$.)  These bounds follow easily from the definitions of $\mu$ and $\nu$; see Section \ref{sec:bounds}.

The operation $\mu$ \cite{Cahn} and Turaev's cobracket $\nu$ \cite{Turaev} were previously defined on the free $\Z$-module generated by the set of nontrivial free homotopy classes of loops on a surface.  In \cite{Cahn}, we proved that, when formulated for curves on an orientable surface $F$, Inequality \eqref{virtualmubound} is an equality.  In other words, we have
$$m(\langle\alpha\rangle)= t(\mu(\langle\alpha\rangle))/2 +n-1,$$
where $m(\langle \alpha\rangle)$ denotes the minimum number of self-intersection points of a curve in the free homotopy class $\langle \alpha \rangle$, and $n\geq 1$ is the largest integer such that $\langle \alpha \rangle = \langle \gamma \rangle ^n$ for some $\langle \gamma \rangle \in \pi_1(F)$.  On the other hand, Chas \cite{Chas} showed that, when formulated for curves on surfaces, Inequality \eqref{virtualnubound} is not an equality in general.  

We ask whether the above bounds are equalities in the virtual category.  We suspect that $m([\alpha])= t(\mu([\alpha]))/2 +n-1$ is true in the virtual category for all classes $[\alpha]$; Theorem \ref{primitiveformulaintro} allows us to construct examples of classes $[\alpha]$ such that this equality holds.  On the other hand, there are strings such that $m([\alpha])\neq t(\nu([\alpha]))/2 +n-1$ in the virtual category; see Section \ref{sec:bounds} and Corollary \ref{primitiveformulaintrocor}.   

Turaev's primitive based matrices \cite{Turaev2}, which appear in the statement of Theorem \ref{primitiveformulaintro}, are reviewed in Section \ref{sec:matrices}.  Self-complementary elements are defined in Subsection \ref{sec:compositemoves}.  \\\\

\begin{thm} \label{primitiveformulaintro} Let $\alpha$ be a virtual string, whose based matrix $T(\alpha)$ is primitive and does not contain a self-complementary element.  Then $m([\alpha])=t(\mu([\alpha]))/2$.    If $T(\alpha)$ is primitive and does contain a self-complementary element, then either $m([\alpha])=t(\mu([\alpha]))/2$ or $m([\alpha])=t(\mu([\alpha]))/2+1$ .
\end{thm}

 (Note that for the strings where $m([\alpha])=t(\mu([\alpha]))/2$, we must have $n=0$, and for strings where $m([\alpha])=t(\mu([\alpha]))/2+1$, we must have $n=0$ or $n=1$.) 

It follows from this result that, as is the case for curves on surfaces, the bound on the minimal self-intersection number given by $\mu$ is stronger than the bound given by $\nu$.

\begin{cor} \label{primitiveformulaintrocor}  The bound on $m([\alpha])$ given by $\mu$ is stronger than the bound given by Turaev's cobracket $\nu$.  Namely, the number of terms of $\mu([\alpha])$ is greater than or equal to the number of terms of $\nu[\alpha]$, and there are virtual homotopy classes $[\alpha]$ such that this inequality is strict.
\end{cor}

Theorem \ref{main2} compares the bound on $m([\alpha])$ given by $\mu$ to a lower bound $\rho$ given by Turaev's based matrix invariant.  The integer-valued invariant $\rho$ is introduced in \cite{Turaev2}, and we discuss it in Section \ref{sec:matrices}.  For now we just note that $m([\alpha])\geq \rho([\alpha])$.  The quantity $O$ is an even integer that counts certain elements of the based matrix of $\alpha$ and is defined in Section \ref{sec.boundcomparison}.

\begin{thm} \label{main2} Let $\alpha$ be any virtual string.  Then 
	$$m([\alpha])\geq t(\mu(\alpha))/2\geq \rho([\alpha])-1+O/2.$$
	If the primitive based matrix associated to $[\alpha]$ does not contain a self-complementary element, then 
	$$m([\alpha])\geq t(\mu(\alpha))/2\geq \rho([\alpha])+O/2.$$
	
\end{thm}

While we can only prove that $t(\mu([\alpha]))/2\geq\rho([\alpha])-1$ in general, and hence that $t(\mu([\alpha]))/2+n-1\geq\rho([\alpha])-1$, we do not know of any examples where $t(\mu([\alpha]))/2+n-1=\rho([\alpha])-1$, and we suspect that $t(\mu([\alpha]))/2+n-1\geq\rho([\alpha])$.

In fact, we construct a string $\alpha$ that shows that the bound given by $\mu$ is sometimes stronger than the bound $\rho([\alpha])$.  The same string also illustrates that the bound given by $\mu$ is stronger than the bound given by $\nu$.
\begin{mex} \label{MainExample}The string $\alpha$ in Figure \ref{Gibsonex.fig} satisfies the following inequality: $t(\mu([\alpha]))/2 > \rho([\alpha])$.  More precisely, $t(\mu([\alpha]))/2=5$ and $\rho([\alpha])=4$.  It is also easy to check that $\nu([\alpha])=0$, so $t(\mu([\alpha]))/2>t(\nu([\alpha]))/2$, illustrating Corollary \ref{primitiveformulaintrocor}.
	\begin{figure}[h]
	   \centering
	   \includegraphics[width=6cm]{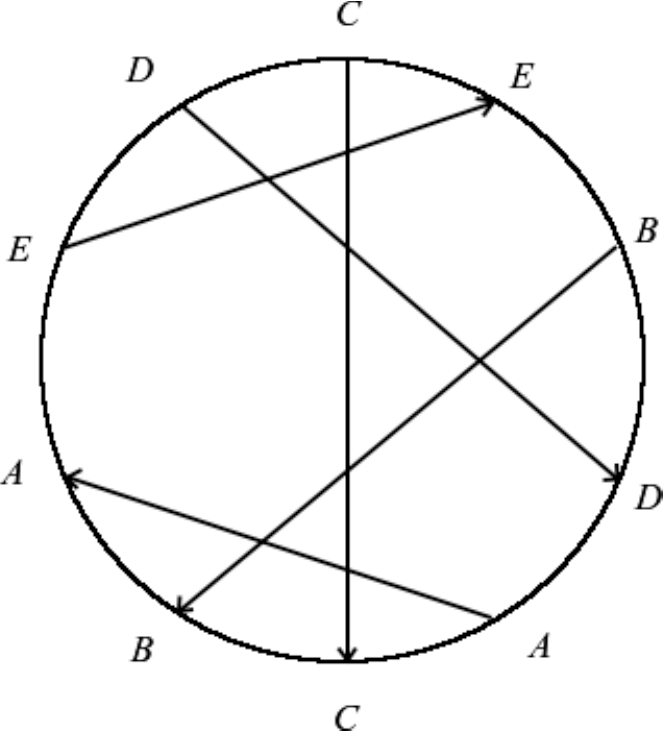} 
	  \caption{A string $\alpha$ such that the bound on $m([\alpha])$ given by $\mu([\alpha])$ is greater than $\rho([\alpha])$.}
	   \label{Gibsonex.fig}
	\end{figure}
\end{mex}

{\it Remark:} We do not focus on algorithmic methods of computing $m([\alpha])$ in this paper.   The results of Ilyutko, Manturov, and Nikonov in \cite{IlyutkoManturovNikonov} imply that an algorithm described in \cite{Gibson} compues the minimal self-intersection number.  We use these results to show that the string in Example \ref{MainExample} has the desired properties.

\section{Virtual Strings, Singular Virtual Strings, and Signed Singular Virtual Strings}\label{sec:strings}
In this section, we recall Turaev's definition of virtual strings and their homotopy classes, and we use the notation found in \cite{Turaev2}.  We then define signed singular virtual strings and their homotopy classes, which are closely related to the singular virtual strings studied by Henrich \cite{Henrich}.  Once we have these definitions, we will be able to define Turaev's cobracket and the operation $\mu$.  Turaev's cobracket takes values in the free $\Z$-module generated by nontrivial homotopy classes of virtual strings, tensored with itself over $\Z$.  The operation $\mu$ takes values in the free $\Z$-module generated by homotopy classes of signed singular virtual strings.

\begin{defin}[Turaev \cite{Turaev2}]  A {\it virtual string $\alpha$ of rank $m$} is an oriented copy of $S^1$, called the core circle of $\alpha$, with $m$ arrows $(a,b)$ whose tail and head are attached to the points $a$ and $b$ of $S^1$, respectively.  We require that the endpoints of the arrows be distinct.  
\end{defin}
\noindent The set of arrows of $\alpha$ is denoted $\arr(\alpha)$.  

\begin{defin}[Henrich, p. 5  \cite{Henrich}] A {\it singular virtual string} $\alpha_d$ is a virtual string with a choice of distinguished arrow $d\in \arr(\alpha)$.
\end{defin}

\begin{defin}A {\it signed singular virtual string} $\alpha_d^\epsilon$ is a singular virutal string with a distinguished arrow $d$, where the distinguished arrow is equipped with a sign $\epsilon \in \{+,-\}$.
\end{defin}

\subsection{The underlying (singular) string of a (singular) curve.}  
A {\it curve} on an oriented surface $F$ is a generic immersion $A: S^1\rightarrow F$, where generic means that its self-intersection points are transverse double points.  Every curve on an oriented surface has an {\it underlying virtual string}, which is also known as its Gauss diagram.  This string $u(A)$ is obtained as follows.  There is one arrow $(a,b)$ corresponding to each self-intersection point $p$ of $A$. The tail $a$ and head $b$ are chosen so that the ordered pair of tangent vectors $\{A'(a),A'(b)\}$ gives a positive orientation of $T_pF$. The string $u(A)$ depends on the choice of orientation of $F$; changing the orientation of $F$ would reverse the direction of each arrow in $u(A)$.

We say that the curve $A$ {\it realizes} the string $\alpha$ if $\alpha$ is the underlying string of $A$, i.e., $\alpha=u(A)$. Conversely, every string $\alpha$ realizes some curve $A$ on an oriented surface $F$, and this can be done canonically; this construction can be found in \cite[p. 2468]{Turaev2}.

Let $\Theta$ be a copy of $S^1$ with an interval $I=[0,1]$ (called a chord) attached to $S^1$ at its endpoints.  A {\it singular curve} on an oriented surface $F$ is a map $A:\Theta\rightarrow F$ such that $A$ maps the chord of $\Theta$ to a single point, and such that $A|_{S^1}$ is a curve on $F$.  When we draw $A$, we draw a thick dot on top of the self-intersection point of $A$ whose preimage is the chord of $\Theta$, and we call this the {\it distinguished} self-intersection point of $A$.  A singular curve is a special case of a {\it geometrical chord diagram}, defined by Andersen, Mattes, and Reshetikhin \cite{AMR}.

A {\it signed singular curve} on an oriented surface $F$ is a singular curve $A$ on $F$ whose distinguished self-intersection point is equipped with a sign $\epsilon \in \{+,-\}$.  

The {\it underlying (signed) singular virtual string} $u(A)$ of a (signed) singular curve $A$ is the underlying virtual string of $A|_{S^1}$, whose distinguished arrow is the arrow corresponding to the distinguished self-intersection point of $A$ (and is equipped with the sign of the distinguished self-intersection point of $A$ if $A$ is signed).

\subsection{Homotopies of signed singular curves.}

Two curves on $F$ are in the same free homotopy class if and only if they are related by a sequence of flat Reidemeister moves.

\begin{figure}[h]
   \centering
   \includegraphics[width=3cm]{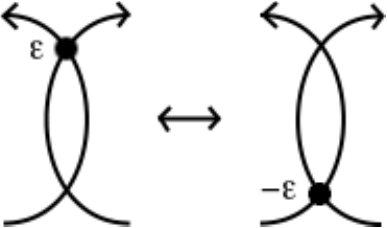} 
  \caption{Flat Reidemeister move for a signed singular curve.  There is a version of this move for each choice of orientation on the two arcs.}
   \label{Type2movewithdot.fig}
\end{figure}

\begin{figure}[h]
	   \centering
	   \includegraphics[width=4cm]{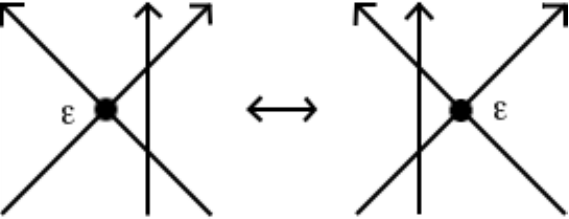} 
	  \caption{Flat Reidemeister move for a signed singular curve.  There is a version of this move for each choice of orientation on the three arcs.}
	   \label{Type3movewithdot.fig}
	\end{figure}

We say two signed singular strings are homotopic if and only if they are related by flat Reidemeister moves which do {\bf not} involve the distinguished self-intersection point, along with the moves in Figures \ref{Type2movewithdot.fig} and \ref{Type3movewithdot.fig}.  

	

Every signed singular curve $A$ can be regarded as a map $\bar{A}:S^1\vee S^1\rightarrow F$, such that the two copies of $S^1$ are ordered; i.e., labeled with a `1' and a `2'.  It is clear how to view $A$ as a map from an unordered wedge of two circles to $F$, since if we collapse the chord of $\Theta$ to a point, we get $S^1\vee S^1$.  It remains to explain how to order the two copies of $S^1$.  Since $F$ is oriented and $A$ is generic, we order the pair tangent vectors to $A$ at the distinguished self-intersection point $p$ of $A$ according to the orientation of $F$.  This gives us an order on the two outgoing branches of $A$ at $p$.  If the sign $\epsilon$ of $p$ is $+$ (respectively, $-$), we label the loop starting on the first (respectively second) outgoing branch with a `1', and the loop starting on the second (respectively first) outgoing branch with a `2'.  Let $\bar{A}_i\in \pi_1(F)$ denote the restriction of $\bar{A}$ to the loop labeled $`i$', $i=1,2$.

Later, we will use the following proposition, whose proof is straightforward:
\begin{prop}\label{figure8s} Two signed singular curves $A$ and $B$ are homotopic as signed singular curves if and only if $\bar{A}$ and $\bar{B}$ are free homotopic as maps from an ordered wedge of circles $S^1\vee S^1$ to $F$; i.e., if and only if there exists $\gamma\in\pi_1(F)$ such that $\gamma \bar{A}_i \gamma^{-1}=\bar{B}_i$ for $i=1,2$.
\end{prop}

\subsection{The homotopy class of virtual string.}  By performing a flat Reidemeister move on a curve, and recording the effect on its underlying virtual string, one can derive the definition of a Reidemeister move for a virtual strings.  This allows one to define a homotopy of virtual strings.  

Two strings are in the same virtual homotopy class if they can be related by a sequence of the following moves and their inverses \cite[p. 2459]{Turaev2}.  We use Turaev's notation for consistency.\\\\
Type 1:  Given an arc $ab$ of $S^1$ containing no endpoints of $\alpha$, add an arrow $e=(a,b)$ or an arrow $e=(b,a)$.\\\\
Type 2:  Let $a$ and $a'$, and $b$ and $b'$, be the endpoints of two disjoint arcs of $S^1$, such that neither arc contains an endpoint of an arrow of $\alpha$.  Add two arrows $e=(a,b)$ and $e'=(b',a')$ to $\alpha$.  There are four forms of this move depending on the order in which $a$ and $a'$ appear, and the order in which $b$ and $b'$ appear, as one traverses $S^1$ counterclockwise.  Note that the new arrows point in opposite directions in all forms of this move. \\\\
Type 3a: Let $aa^+$, $bb^+$, and $cc^+$ be three disjoint arcs of $S^1$, containing no endpoints of $\alpha$.    Suppose also that $(a^+,b)$, $(b^+,c)$, and $(c^+, a)$ are arrows of $\alpha$.  Replace these arrows with the arrows $(a,b^+)$, $(b,c^+)$, and $(c,a^+)$.\\\\
Type 3b: Let $aa^+$, $bb^+$, and $cc^+$ be three disjoint arcs of $S^1$, containing no endpoints of $\alpha$. Suppose also that  $(a,b)$, $(a^+,c)$ and $(b^+,c^+)$ are three arrows of $\alpha$.  Replace these arrows with the arrows  $(a^+,b^+)$, $(a,c^+)$, and $(b,c)$.  Turaev notes this move is not necessary because it can be expressed as a composition of a Type 2 and Type 3a move.\\\\
The virtual homotopy class of $\alpha$ is denoted by $[\alpha]$.  

 One can view a virtual string as a curve on a surface up to free homotopy and stable homeomorphism \cite{Turaev2} (the case for ordinary virtual knots appeared in \cite{CarterKamadaSaito}).\\\\

\subsection{The homotopy class of a signed singular virtual string.}
Two singular virtual strings are homotopic if they can be related by a sequence of Type 1-3 moves for virtual strings, that do {\bf not} affect the distinguished arrow, along with the following moves, and their inverses:\\\\
Signed Singular Type 2:  Let $a$ and $a'$, and $b$ and $b'$, be the endpoints of two disjoint arcs of $S^1$ that contain no endpoints of arrows of $\alpha$.   Suppose $e=(a,b)$ and $e'=(b',a')$ are arrows of $\alpha$, and suppose furthermore that one of these is the distinguished arrow.  This move changes which of $e$ and $e'$ is the distinguished arrow, and changes the sign of the distinguished arrow.  Note there are four forms of this move.  Also note that no arrows are removed (in contrast to the ordinary Type 2 move).

Signed Singular Type 3a and 3b:  This is the same as the Type 3a or 3b move for ordinary virtual strings, but where one of the three arrows, and the arrow which replaces it, is the distinguished arrow, and the sign on the distinguished arrow is the same before and after the move.\\\\

The virtual homotopy class of a signed singular string $\alpha^\epsilon_d$ is denoted $[\alpha^\epsilon_d]$.

If one forgets the signs on these moves, one recovers the moves for singular virtual strings defined in \cite{Henrich}.

It is easy to see that the underlying signed singular virtual strings of two homotopic signed singular curves on $F$ are virtually homotopic, and that the underlying virtual strings of two homotopic curves on $F$ are virtually homotopic.

\subsection{Semi-trivial signed singular strings.} We call a signed singular virtual string $\textit{semi-trivial}$ if the endpoints $a$ and $b$ of the distinguished arrow form an arc $ab$ or $ba$ of $S^1$ containing no endpoints in its interior.  We call a homotopy class of a signed singular virtual string $\textit{semi-trivial}$ if it contains a semi-trivial signed singular string.\\\\

\section{Turaev's Cobracket for Virtual Strings} 
\noindent Let $e=(a,b)$ be an arrow of the virtual string $\alpha$.  Let $\alpha^1_e$ (respectively, $\alpha^2_e$) be the virtual string obtained from $\alpha$ by deleting all arrows except those with tail and head in the interior of the oriented arc $ab$ (respectively, $ba$). Let $\mathcal{S}_0$ be the set of nontrivial homotopy classes of virtual strings.  Turaev's cobracket $\nu: \Z[\mathcal{S}_0]\rightarrow\Z[\mathcal{S}_0]\otimes\Z[\mathcal{S}_0]$ is a linear map, defined on a single class by
$$\nu([\alpha])=\sum_{e\in \arr(\alpha)}[\alpha^1_e]\otimes[\alpha^2_e]-[\alpha^2_e]\otimes[\alpha^1_e],$$
where we set $[\beta]=0$ if $\beta$ is trivial.  This map can be extended by linearity to all of $\Z[\mathcal{S}_0]$.  One can verify that this definition is independent of the choice of $\alpha\in [\alpha]$.  Furthermore, the operation $\nu$ is a Lie cobracket \cite{Turaev2}.
\section{The Operation $\mu$}
Let $\mathcal{S}^\Theta_0$ denote the set of homotopy classes of singular virtual strings which are not semi-trivial.  Set $[\beta]=0$ if $\beta$ is semi-trivial. 
\noindent We define 
$$\mu([\alpha])=\sum_{e\in \arr(\alpha)} [\alpha_e^+]-[\alpha_e^-].$$
  We extend the definition of $\mu$ by linearity to obtain a map $\mu: \mathbb{Z}[\mathcal{S}_0]\rightarrow \mathbb{Z}[\mathcal{S}^\Theta_0]$.
\subsection{$\mu$ is well-defined}  
\begin{figure}[htbp]
   \centering
   \includegraphics[width=8cm]{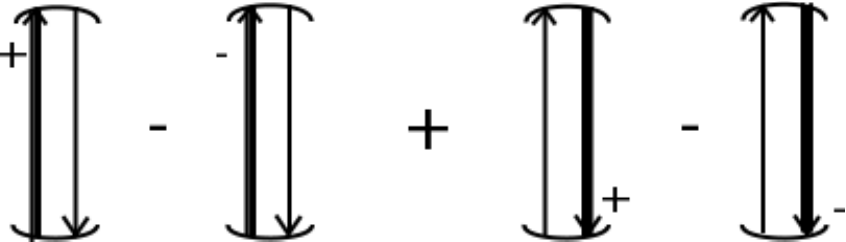} 
  \caption{The signed singular Type 2 move allows one to prove that $\mu$ is invariant under the Type 2 move.}
   \label{muinvarR2.fig}
\end{figure}

Suppose we compute $\mu$ before and after the Type 1 move.  After the move, the signed singular virtual strings which come from the new arrow will be semi-trivial, so they do not contribute anything to the sum.\\\\ 
Figure \ref{muinvarR2.fig} shows the four terms which one form of the Type 2 move contributes to $\mu$. The first and last terms cancel after an application of the signed singular Type 2 move, as do the second and third terms.  The arguments for the other forms of this move are similar.

Similarly the singular Type 3 move allows us to show that $\mu$ is invariant under the ordinary Type 3 move.
\subsection{$\nu$ Factors through $\mu$}  \label{sec.smoothingmap}Let $S$ be the following ``smoothing" map: \\\\
Given a signed singular virtual string $\alpha$ with distinguished arrow $d=(t,h)$ and sign $\epsilon=+$ (respectively $\epsilon=-$), form the string $S_1(\alpha)$ by deleting the distinguished arrow as well as all arrows except those with head and tail in $(th)^\circ$ (respectively $(ht)^\circ$), and form the string $S_2(\alpha)$ by deleting the distinguished arrow as well as all arrows except those with head and tail in $(ht)^\circ$ (respectively $(th)^\circ$).  Let $S(\alpha)=S_1(\alpha)\otimes S_2(\alpha)$, and extend $S$ by linearity to $\Z[\mathcal{S}^\Theta_0]$.  Then
$$\nu(\beta)=S\circ \mu(\beta).$$
 \begin{figure}[h]\center
\scalebox{1.0}{\includegraphics[width=4cm]{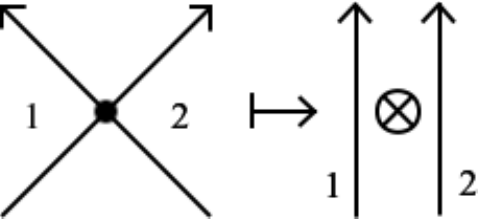}}  
\caption{Illustration of the smoothing map $S(\alpha)$ after realizing the signed singular string $\alpha$ on a surface.}
\label{factors1.fig}
\end{figure}

\section{Bounds on the Minimal Self-Intersection Number Given by $\mu$ and $\nu$}\label{sec:bounds}
\begin{figure}[h]\center
\scalebox{1.0}{\includegraphics[width=4cm]{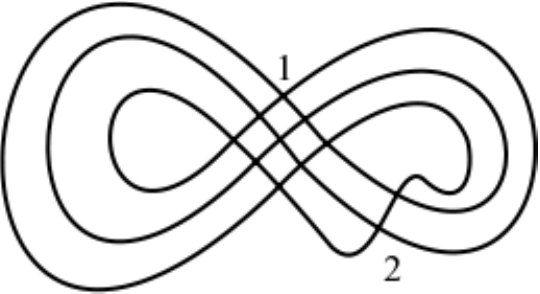}}  
\caption{The composition $h=\bar{g}\circ d_n$ where $n=3$ and $g$ is a figure-8 in the plane.}
\label{power.fig}
\end{figure}

In this section we show why $\mu$ and $\nu$ give the bounds on $m([\alpha])$ stated in the introduction.
\begin{prop} \label{powers} Let $\alpha$ be a virtual string.  Then
$$m([\alpha])\geq t(\mu([\alpha]))/2 +n-1$$ and
$$m([\alpha])\geq t(\nu([\alpha]))/2 +n-1,$$ where $n\geq 1$ is the largest integer such that, for some $\beta\in[\alpha]$ with minimal self-intersection, $\beta$ can be realized as a curve $B$ on an orientable surface $F$ with $\langle B \rangle=\langle \gamma \rangle ^n$ for some $\langle \gamma \rangle$ in $\pi_1(F)$. 
\end{prop} 
\pp
Suppose that $n\geq 1$ is the largest integer such that, for some $\beta\in[\alpha]$ with minimal self-intersection, $\beta$ can be realized as a curve $B$ on a surface $F$ with $\langle B \rangle=\langle \gamma \rangle ^n$ for some $\langle \gamma \rangle$ in $\pi_1(F)$.  We will prove that $m([\alpha])\geq t(\mu([\alpha]))/2 +n-1$; the argument for $\nu$ is slightly simpler.  Let $g$ be a geodesic representative of $\langle \gamma \rangle$.  If $g$ is not generic (i.e., if it has self-intersection points which are not double points), perturb $g$ slightly to get a generic loop $g'$ and put $g=g'$.  Let $h$ be the composition of the following maps: $d_n:S^1\rightarrow A$, the degree $n$ map of $S^1$ to the annulus $A$ such that the image of $d_n$ has $n-1$ double points; and $\bar{g}:A\rightarrow F$, an immersion of the annulus into a thin neighborhood of the loop $g$.  See Figure \ref{power.fig} for an illustration of this composition when $n=3$.  Hass and Scott \cite[p. 10]{HassScott} show this loop $h$ has the fewest number of self intersection points of any loop in its free homotopy class on $F$.  Since $\beta$ has minimal self-intersection in its virtual homotopy class, the loop $B$ on $F$ realizing $\beta$ has minimal self-intersection in its free homotopy class.  So the number of self intersection points of $B$ is equal to the number of self intersection points of $h$, and hence we may assume $\beta$ is the underlying virtual string of $h$, i.e., we assume $B=h$.  The number of self intersection points of $h$ is $n^2p+n-1$, where $p$ is the number of self-intersection points of $g$.  

Now we compute $\mu([\alpha])$ using the representative $\beta=u(h)$.  There are two terms of $\mu$ for each self intersection point $d$ of $h$.  These two terms are $+[\beta^+_d]$ and $-[\beta^-_d]$.  The signed singular virtual strings $\beta^+_d$ and $\beta^-_d$ are realized by the signed singular virtual curves $h^+_d$ and $h^-_d$ on $F$.  Now suppose $d$ is one of the $n-1$ self-intersection points of $h$ that comes from the map $d_n:S^1\rightarrow A$; these self-intersection points are labeled with a `2' in Figure \ref{power.fig}. Starting on the inside of the annulus, we label these $n-1$ points $e_1,\dots, e_{n-1}$.  It is straightforward to check that the signed singular curve $h^+_{e_i}$ is free homotopic to the singular virtual string $h^-_{e_{n-i}}$, and the singular curve $h^-_{e_i}$ is free homotopic to the singular virtual curve $h^+_{e_{n-i}}$, using Proposition \ref{figure8s}.  Therefore their underlying signed singular virtual strings are virtually homotopic.  Thus the $2(n-1)$ terms corresponding to the self intersection points $e_1,\dots e_{n-1}$ cancel. Now $n^2p\geq t(\mu([\alpha]))/2$, so $m([\alpha])=n^2p+n-1\geq t(\mu([\alpha]))/2+n-1$.

\qed
\section{Signed Singular Based Matrices}\label{sec:matrices} In this section, we introduce signed singular based matrices.  These will be our main tools in the proofs of Theorems 1.1 and 1.3.   Signed singular based matrices are closely related to Turaev's based matrices \cite{Turaev2} and Henrich's singular based matrices \cite{Henrich}.

\subsection{Based matrices}  First we recall Turaev's definition of a based matrix.  A \textit{based matrix} over an abelian group $H$ is a triple $(G,s,b)$ where $G$ is a finite set, $s\in G$, and $b:G\times G\rightarrow H$ is a skew-symmetric map.  That is, $b(g,h)=-b(h,g)$ for all $g,h\in G$ and $b(g,g)=0$ for all $g\in G$.  Certain elements of $G$ allow us to perform moves on based matrices that are analogous to flat Reidemeister moves.  These elements are:\\
1) Annihilating elements: An element $g\in G-\{s\}$ is \textit{annihilating} if $b(g,h)=0$ for all $h\in G$.\\
2) Core elements: An element $g\in G-\{s\}$ is \textit{core} if $b(g,h)=b(s,h)$ for all $h\in G$.\\
3) Complementary elements: Two elements $g_1, g_2\in G-\{s\}$ are \textit{complementary} if $b(g_1,h)+b(g_2,h)=b(s,h)$ for all $h\in G$.\\

\subsection{Signed singular based matrices}  A signed singular based matrix is a quadruple $(G,s,d,b,\epsilon)$ where $(G,s,b)$ is a based matrix, $d\in G-\{s\}$, and $\epsilon\in \{+,-\}$.  In addition to annihilating, core, and complementary elements, all of which are required to be elements other than $s$ and $d$, we also have the following:\\
1)  Annihilating-like elements: The distinguished element $d$ is \textit{annihilating-like} if $b(d,h)=0$ for all $h\in G$.  Similarly, $s$ is annihilating-like if $b(s,h)=0$ for all $h\in G$.\\
2) Core-like elements: The distinguished element $d$ is \textit{core-like} if $b(d,h)=b(s.h)$ for all $h\in G$.

If one forgets the sign in the definitions above, one recovers Henrich's singular based matrices \cite{Henrich}.

\subsection{Elementary extensions of based matrices.}  Turaev defined the following moves, called {\it elementary extensions}, on based matrices \cite[p. 2483]{Turaev2}:\\\\
$M_1$ adds an annihilating element to $(G,s,b)$:  Given $(G,s,b)$, form the based matrix $(\bar{G},s,\bar{b})$ where $\bar{G}=G\sqcup\{g\}$, $\bar{b}$ agrees with $b$ on $G\times G$, and $\bar{b}(g,h)=0$ for all $h\in G$.\\\\
$M_2$ adds a core element to $(G,s,b)$:  Given $(G,s,b)$, form the based matrix $(\bar{G},s,\bar{b})$ where $\bar{G}=g\sqcup\{g\}$, $\bar{b}$ agrees with $b$ on $G\times G$, and $\bar{b}(g,h)=\bar{b}(s,h)$ for all $h\in G$.\\\\
$M_3$ adds a pair of complementary elements to $(G,s,b)$:  Given $(G,s,b)$, form the based matrix $(\bar{G},s,\bar{b})$, where $\bar{G}=G\sqcup \{g_1,g_2\}$, $\bar{b}$ agrees with $b$ on $G$, and $\bar{b}(g_1,h)+\bar{b}(g_2,h)=\bar{b}(s,h)$ for all $h\in G$.\\\\
The inverses of these moves are called {\it inverse extensions}.
\subsection{Elementary extensions of signed singular based matrices.}
For signed singular based matrices, we define the following elementary extensions, plus an additional move which changes which element is the distinguished element and changes the sign of the matrix:\\\\
$M''_1$ adds an annihilating element to $(G,s,d,b, \epsilon)$:  Given $(G,s,d,b,\epsilon)$, form the signed singular based matrix $(\bar{G},s,d,\bar{b},\epsilon)$ where $\bar{G}=G\sqcup\{g\}$, $\bar{b}$ agrees with $b$ on $G\times G$, and $\bar{b}(g,h)=0$ for all $h\in G$.\\\\
$M''_2$ adds a core element to $(G,s,d,b,\epsilon)$:  Given $(G,s,d,b,\epsilon)$, form the signed singular based matrix $(\bar{G},s,d,\bar{b},\epsilon)$ where $\bar{G}=g\sqcup\{g\}$, $\bar{b}$ agrees with $b$ on $G\times G$, and $\bar{b}(g,h)=\bar{b}(s,h)$ for all $h\in G$.\\\\
$M''_3$ adds a pair of complementary elements to $(G,s,d,b,\epsilon)$:  Given $(G,s,d,b,\epsilon)$, form the singular based matrix $(\bar{G},s,d,\bar{b},\epsilon)$, where $\bar{G}=G\sqcup \{g_1,g_2\}$, $\bar{b}$ agrees with $b$ on $G$, and $\bar{b}(g_1,h)+\bar{b}(g_2,h)=\bar{b}(s,h)$ for all $h\in G$.\\\\
$N''$ changes which element is the distinguished element, as well as the sign of the matrix:  Given $(G,s,d,b,\epsilon)$ such that $g\in G$ and $d$ are complementary, form the singular based matrix $(G,s,g,b,-\epsilon)$.\\\\
By forgetting the signs on the matrices, one obtains the moves $M'_i$ and $N'$ for singular based matrices in \cite{Henrich}.\\\\

\subsection{Homologous, primitive, and isomorphic matrices} 
Two based matrices are \textit{homologous} if one can be obtained from the other by a finite number of the $M_i$ moves and their inverses \cite{Turaev2}.  Two signed singular based matrices are homologous if one can be obtained from the other by a finite number of $M''_i$ moves, $N''$ moves, and their inverses.  

A based matrix is \textit{primitive} if it cannot be obtained from another based matrix by a sequence of $M_i$ moves \cite{Turaev2}.  A signed singular based matrix is {\it primitive} if it cannot be obtained from another singular based matrix by applying an $M''_i$ move, possibly preceded by an $N''$ move.

Two based matrices $(G,s,b)$ and $(G',s',b')$ are \textit{isomorphic} if there is a bijection $\phi: G\rightarrow G'$ such that $\phi(s)=s'$, and $\phi(b(g,h))=b'(\phi(g),\phi(h))$ for all $g,h\in G$. Two signed singular based matrices $(G,s,d,b,\epsilon)$ and $(G',s',d',b',\epsilon)$ are \textit{isomorphic} if there is a bijection $\phi: G\rightarrow G'$ such that $\phi(s)=s'$, $\phi(d)=d'$, $\epsilon=\epsilon'$, and $\phi(b(g,h))=b'(\phi(g),\phi(h))$ for all $g,h\in G$.

By forgetting signs, one recovers the definitions of homologous, primitive, and isomorphic singular based matrices in \cite{Henrich}.
\section{Associating a Matrix to a Virtual String}
We will recall how to associate a based matrix to a virtual string \cite{Turaev2}, and define a method of associating a signed singular based matrix to a signed singular virtual string.  By forgetting signs, one recovers the definitions of \cite{Henrich}.
\subsection{Associating a based matrix to a virtual string}  
Turaev gives a combinatorial formula for the based matrix of a virtual string \cite{Turaev2}.

\begin{figure}[htbp]
   \centering
   \includegraphics[width=5cm]{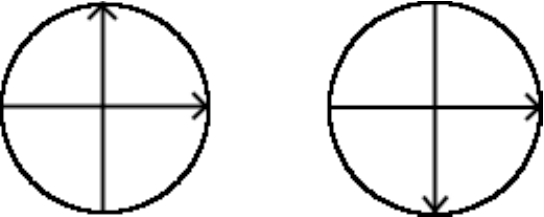} 
  \caption{In both figures $e$ is the horizontal arrow and $f$ is the vertical arrow.  On the right, $f$ links $e$ positively.  On the left, $f$ links $e$ negatively.}
   \label{linkingdefn.fig}
\end{figure}

The based matrix $T(\alpha)=(G, s, b)$ of the string $\alpha$ is given as follows: Put $G=\{s\}\cup arr(\alpha)$.
Suppose $e=(a,b)$ and $f=(c,d)$ are arrows of $\alpha$.  We say $f$ links  $e$ \textit{positively} if $c$ lies in the arc $ab$ and $d$ lies in the arc $ba$.  We say $f$ links $e$ \textit{negatively} if $c$ lies in the arc $ba$ and $d$ lies in the arc $ab$.  Otherwise, $e$ and $f$ are \textit{unlinked}. (See Figure \ref{linkingdefn.fig}.)
 Let
$$n(e)=\#\{f\in \arr(\alpha)|f \text{ links } e \text{ positively}\}-\#\{f\in \arr(\alpha)|f \text{ links } e \text{ negatively}\}.$$
Then put $b(e,s)=n(e)$.  (This defines the first row and column of the matrix.) Now let $(ab)^\circ$ denote the interior of the arc $ab$, and for two arcs $ab$ and $cd$, and let 
$$ab\cdot cd=\#\{(t,h)\in\arr(\alpha)|t\in (ab)^\circ \text{ and } h\in (cd)^\circ\}$$
$$-\#\{(t,h)\in\arr(\alpha)|t\in (cd)^\circ \text{ and } h\in (ab)^\circ\}.$$ 
Set $\sigma=0$ if $e$ and $f$ are unlinked, set $\sigma=1$ if $f$ links $e$ positively, and set $\sigma=-1$ if $f$ links $e$ negatively.  Then put $b(e,f)=ab\cdot cd +\sigma$.\\\\
\subsection{Associating a signed singular matrix to a signed singular string.}  Let $\alpha^\epsilon_d$ be a signed singular virtual string, with underlying virtual string $\alpha$, distinguished arrow $d$ and sign $\epsilon$.  The signed singular based matrix $T(\alpha^\epsilon_d)$ associated to $\alpha^\epsilon_d$ is $(G,s, d, b, \epsilon)$, where $T(\alpha)=(G,s,b)$.  When writing $T(\alpha^\epsilon_d)$ in matrix form, we will display the row and column corresponding to $d$ in bold, and display the sign $\epsilon$ to the upper right of the matrix.

\subsection{Example.} \label{sec:matrixofastring} We compute the based matrix of the virtual string $M$ in Figure \ref{Gibsonex.fig}.  Columns one through six correspond to $s,A,B,C,D,$ and $E$ respectively.  We will compute two entries explicitly: $b(A,s)=n(A)=2$, so the matrix entry $a_{21}=2$.  For an entry not in the first row or column, we use the formula $b(e,f)=ab\cdot cd +\sigma$ to get $a_{25}=b(A,D)=2-1+0=1$.  Notice that $C$ is annihilating.  Gibson \cite{Gibson} computes this based matrix, but it seems we have opposite sign conventions, so his matrix is the transpose of ours.  Below, we display the entire matrix $T(M)$, as well as the signed singular matrix of the signed singular string $M^+_C$:
$$T(M)=\left(\begin{matrix}
0 & -2 & -1 & 0 & 1 & 2\\
2 & 0 & 0 & 0 &1 & 3\\
1 & 0 & 0 & 0 & 0 & 1\\
0 & 0 & 0 & 0 & 0 & 0\\
-1 & -1 & 0 & 0 & 0 & 0\\
-2 & -3 & -1 & 0 & 0 & 0\end{matrix}\right), \text{ and }
T(M^+_C)=
\left(\begin{matrix}
0 & -2 & -1 & {\bf 0} & 1 & 2\\
2 & 0 & 0 & {\bf 0} &1 & 3\\
1 & 0 & 0 & {\bf 0} & 0 & 1\\
{\bf 0} & {\bf 0} & {\bf 0} & {\bf 0} & {\bf 0} & {\bf 0}\\
-1 & -1 & 0 & {\bf 0} & 0 & 0\\
-2 & -3 & -1 & {\bf 0} & 0 & 0\end{matrix}\right)^+.$$
\section{Primitive Signed Singular Based Matrices Yield Invariants of Signed Singular Strings}
\subsection{The homology class of a (signed singular) based matrix is an invariant of (signed singular) strings.}In this section we find analogues of Turaev's results for based matrices and Henrich's results for singular based matrices:
\begin{prop} [Turaev, p. 2487] Homotopic virtual strings have homologous based matrices.\label{homotopicimplieshomologous}
\end{prop}
\begin{prop}[Henrich, p. 22] Homotopic singular virtual strings have homologous singular based matrices.
\end{prop}
\begin{prop}\label{homotopicimplieshomologoussignedsingular}  Homotopic signed singular strings have homologous signed singular based matrices.
\end{prop}
\noindent\textit{Proof.}  We need to check that applying the ordinary Type 1-3 moves and the signed singular Type 2 and 3 moves to a signed singular virtual string does not change the homology class of the corresponding signed singular based matrix.  For the ordinary Type 1-3 moves, the proof is identical to Henrich's or Turaev's proof.  Applying a Type 1 move changes the corresponding matrix by adding a core or annihilating element; i.e., by applying the elementary extension $M''_1$ or $M''_2$, depending on the direction of the arrow.  Applying the Type 2 move adds a pair of complementary elements to the corresponding matrix; i.e., we apply the elementary extension $M''_3$.  One can check that the ordinary and signed singular Type 3 moves do not change the corresponding matrix or its sign.  Finally, the signed singular Type 2 move changes the matrix by the move $N''$; i.e., one must change the sign of the matrix, and swap the distinguished element with an element complementary to it.
\qed
\begin{prop} \label{primitive}If the based matrix $T(\alpha)$ is primitive, then the signed singular based matrix $T(\alpha_d^\epsilon)$ is primitive.  
\end{prop}
\pp Clear. \qed

\subsection{The unique primitive representative of a homology class and Turaev's invariant $\rho$.} Our main tool for deciding when terms of $\mu$ cancel is an analogue of the following theorem of Turaev \cite{Turaev2}:

\begin{thm}[Turaev]\label{uniqueprim} There is a unique primitive based matrix in each homology class, up to an isomorphism.
	\end{thm}
Turaev's theorem has two important applications which we will use throughout this paper. First, it gives rise to an invariant of virtual strings.  To see why, recall that by Proposition \ref{homotopicimplieshomologous}, homotopic strings $\alpha$ and $\beta$ have homologous based matrices $T(\alpha)$ and $T(\beta)$.  Each can be reduced to a primitive matrix using inverse elementary extensions.  Therefore if we apply as many inverse extensions as possible to these two matrices, we get the same matrix.  Turaev denotes the unique primitive matrix associated to the homotopy class of $\alpha$ by $T_\bullet([\alpha])$, which we may abbreviate to $T_\bullet(\alpha)$. 

 The second application of this theorem is that it gives a lower bound on the minimal self-intersection number.  Namely, put $\rho([\alpha])=|T_\bullet(\alpha)|-1$.  Then it is clear that $m([\alpha])\geq \rho([\alpha])$.  (This is the invariant $\rho$ mentioned in the introduction.)

\subsection{Example.} The unique primitive matrix associated to the homotopy class of the string $M$ in Figure \ref{Gibsonex.fig} is
$$T_\bullet([M])=\left(\begin{matrix}
0 & -2 & -1  & 1 & 2\\
2 & 0 & 0  &1 & 3\\
1 & 0 & 0  & 0 & 1\\
-1 & -1 &  0 & 0 & 0\\
-2 & -3 & -1 &  0 & 0\end{matrix}\right).$$
We obtained this matrix from $T(M)$ by deleting the center row and column, which corresponded to the annihilating element $C$.  Thus $\rho([M])=5-1=4$.  Later we will see that $m([M])=5$, so the bound $m([\alpha])\geq \rho([\alpha])$ is not an equality in this case.  We will also see that half the number of terms of $\mu([M])$ is in fact 5, so the bound given by $\mu$ is an equality.

\subsection{There is not always a unique primitive matrix in the homology class of a signed singular based matrix.}  Now we seek an analogue of Theorem \ref{uniqueprim} for signed singular based matrices.  There is not a unique primitive representative in the homology class of a signed singular matrix, but there are at most two primitive representatives in each class, and furthermore, if there are two such representatives, we can describe how they are related. 

\subsection{The composite moves $D''_{12}$, $D''_{21}$, and $D''_{33}$.} \label{sec:compositemoves} We will prove that two primitive, homologous, signed singular based matrices either differ by an $N''$ move, or by certain composite moves which we call $D''_{12}$, $D''_{21}$, and $D''_{33}$.  The first two are just signed versions of moves in \cite{Henrich}.  However the move $D''_{33}$ does not appear in the unsigned case.  

Let $D''_{12}$ denote the composition $(M''_1)^{-1}\circ N'' \circ (M''_2)$ in the case where the $N''$ move interacts with both the extension and the inverse extension.  That is, first $M''_2$ adds a core element $c$.  In order for $N''$ to affect this element, it must exchange the current distinguished element $d$ with $c$.  This requires $d$ and $c$ to be complementary, so $d$ must be annihilating.  Then $(M''_1)^{-1}$ removes this annihilating element.  {\it The net effect of $D''_{12}$ is that it replaces an annihilating-like distinguished element with a core-like distinguished element and changes the sign of the matrix.}  

Similarly, we let $D''_{21}$ denote the composition $(M''_2)^{-1}\circ N'' \circ (M''_1)$ in the case where the $N''$ move interacts with both the extension and the inverse extension.  {\it The net effect of $D''_{21}$ is that it replaces a core-like distinguished element with an annihilating-like distinguished element and changes the sign of the matrix.}\\\\

Before introducing the final move, we need the following definition:
\begin{defin}  An element $g\in G$ of the based matrix $(G,s,b)$ is {\it self-complementary} if $2b(g,h)=b(s,h)$ for all $h\in G$.
\end{defin}

Let $D''_{33}$ denote the following special case of the composition $(M''_3)^{-1}\circ N'' \circ M''_3$:  Suppose the distinguished element $d$ is self-complementary, the move $M''_3$ adds a pair of self-complementary elements, and $N''$ switches the distinguished element with one of these two new elements.  Then the $(M'')^{-1}_3$ move removes the old distinguished element and the non-distinguished self-complementary element added by the $M''_3$ move.  {\it The net effect of $D''_{33}$ is that it switches the sign of a matrix whose distinguished element is self-complementary}. The move $D''_{33}$ is equivalent to an isomorphism for unsigned matrices. \\\\

It is helpful to note the following:
\begin{prop}\label{annihilatinglike}A core or annihilating element is self-complementary if and only if $s$ is annihilating-like (in which case core and annihilating elements are the same).
\end{prop}
\pp If $2b(c,h)=b(s,h)$ for all $h\in G$ and $b(c,h)=b(s,h)$, then $b(s,h)=0$.  If $2b(c,h)=b(s,h)$ for all $h\in G$ and $b(c,h)=0$, then $b(s,h)=0$.  The reverse direction is similar.
\qed\\\\

We will use an analogue of the following theorem to understand when certain terms of $\mu$ can cancel.
\begin{thm}[Henrich, p. 19] \label{Allison} Two homologous primitive (unsigned) singular based matrices differ either by an isomorphism or by a composition of an isomorphism with a single $D'_{12}$, $D'_{21}$, or $N'$ move.
\end{thm}
\subsection{Sign switches}  We will need to understand whether the isomorphisms in the statement of Henrich's theorem can change the sign of a signed singular based matrix.  Define an operation $S(T)$ which switches the sign of a signed singular based matrix.  That is, if $T=(G,s,d,b,\epsilon)$, then $S(T)=(G,s,d,b,-\epsilon)$.  Note that $S$ is an isomorphism of the underlying singular based matrices, but $S$ is not an isomorphism of signed singular based matrices.  In general $T$ and $S(T)$ need not be homologous.   \\\\

{\it Remark:} It follows immediately from Theorem \ref{Allison} that two signed singular primitive based matrices differ by an isomorphism, a composition of an isomorphism with a single $D''_{12}$, $D''_{21}$, or $N''$ move, or a composition of one of those moves with a move that changes the sign $\epsilon$ of the matrix.  This is because an isomorphism of unsigned based matrices is either an isomorphism of signed based matrices, or is equivalent to a composition of an isomorphism with a sign-switching move.  But in general, applying a sign switch to a signed based matrix might not produce a matrix homologous to the original matrix.  This is why we need to prove a version of Theorem \ref{Allison} for signed matrices, rather than use Theorem \ref{Allison} directly.  
\begin{thm}\label{primitivesignedtheorem}  Two homologous primitive signed singular based matrices differ by an isomorphism, or by a composition of an isomorphism with a single $D''_{12}$, $D''_{21}$, $D''_{33}$, or $N''$ move.
\end{thm}
Before giving the rather technical proof of Theorem \ref{primitivesignedtheorem}, which we postpone until the next section, we give an example.

As in the case of ordinary based matrices, we let $T_\bullet(\alpha)$ denote a primitive matrix in the class of $T(\alpha)$.  If there is more than one such matrix, the choice of matrix will be specified.

\subsection{Example illustrating Theorem \ref{primitivesignedtheorem}.}  Let us consider a signed singular string whose underlying ordinary string is the string $M$ in Figure \ref{Gibsonex.fig}.  Let $M^+_A$ be the signed singular string with distinguished arrow $A$ and sign $+$.  First we compute the signed singular based matrix associated to $M^+_A$:
$$T(M^+_A)=\left(\begin{matrix}
0 & \bf{-2} & -1 & 0 & 1 & 2\\
\bf{2} & \bf{0} & \bf{0} & \bf{0} &\bf{1} & \bf{3}\\
1 & \bf{0} & 0 & 0 & 0 & 1\\
0 & \bf{0} & 0 & 0 & 0 & 0\\
-1 & \bf{-1} & 0 & 0 & 0 & 0\\
-2 & \bf{-3} & -1 & 0 & 0 & 0\end{matrix}\right)^+.$$
Note that this signed singular matrix is {\bf not} primitive.  A primitive representative of its homology class is:
 $$T_\bullet(M^+_A)=\left(\begin{matrix}
0 & \bf{-2} & -1 &  1 & 2\\
\bf{2} & \bf{0}  & \bf{0} &\bf{1} & \bf{3}\\
1 & \bf{0} & 0  & 0 & 1\\
-1 & \bf{-1} & 0 & 0 & 0\\
-2 & \bf{-3} & -1 &  0 & 0\end{matrix}\right)^+.$$
This matrix happens to be the {\bf unique} primitive representative of its homology class.  Indeed, by Theorem \ref{primitivesignedtheorem}, any other primitive signed singular matrix homologous to this matrix (which is not isomorphic to it) would differ from this one by a $D_{12}''$, $D_{21}''$, $N''$, or $D_{33}''$ move.  But in order to apply these moves, the distinguished element must be either core, annihilating, part of a complementary pair, or self-complementary, and it is easy to check that this is not the case. 

On the other hand, the matrix $T(M^+_C)$ computed in Example \ref{sec:matrixofastring} is primitive (unlike $T(M^+_A)$), and is {\bf not} the unique primitive matrix in its class because its distinguished element is annihilating.  The other primitive matrix in the class of $T(M^+_C)$ is $D''_{12}(T(M^+_C))$.

\subsection{The primitive matrix in a signed singular homology class is sometimes unique.} \label{sec.unique} This example illustrates that there is sometimes a unique primitive matrix in the homology class of a signed singular matrix.  We summarize this in Corollary \ref{uniquematrix} of Theorem \ref{primitivesignedtheorem} below.

Let $\alpha$ be a virtual string.  There may be more than one possible sequence of inverse extensions which reduces the based matrix $T(\alpha)$ to the primitive based matrix $T_\bullet(\alpha)$.  We fix one such sequence $R$, and we let $P(R)$ be the set of arrows of $\alpha$ which are not removed by the reduction $R$.  The set $P(R)$ forms a primitive submatrix of $T(\alpha)$, and we can identify its elements with those of $T_\bullet(\alpha)$. 
\begin{cor}\label{uniquematrix}  For any fixed reduction $R$ of $T(\alpha)$, and any $e\in P(R)$, one can construct a primitive signed singular based matrix $T_\bullet(\alpha^\epsilon_e)$ in the homology class of $T(\alpha^\epsilon_e)$ by using $T(\alpha)$ as the underlying ordinary based matrix, with distinguished element $e$ and sign $\epsilon$.  Furthermore, if $e$ is not a self-complementary element of $P(R)$, then  $T_\bullet(\alpha^\epsilon_e)$ is the unique primitive signed singular based matrix in its homology class.
\end{cor}
\pp  Form the matrix $T(\alpha^\epsilon_e)$.  Since $e\in P(R)$, no move in $R$ removes $e$.  Thus we can apply the entire sequence of moves $R$ to $T(\alpha^\epsilon_e)$ to obtain a matrix we call $T_\bullet(\alpha^\epsilon_e)$.  The signed singular based matrix $T_\bullet(\alpha^\epsilon_e)$ has underlying ordinary based matrix $T(\alpha)$, distinguished element $e$, and sign $\epsilon$.  By Proposition \ref{primitive}, this signed singular based matrix is primitive.  Because $e\in P(R)$, $e$ is not core, annihilating, or part of a complementary pair.   Thus one cannot apply the moves $D''_{12}$, $D''_{21}$, or $N''$ to $T_\bullet(\alpha^\epsilon_e)$.  Hence if $e$ is not a self-complementary element of $P(R)$, then $T_\bullet(\alpha^\epsilon_e)$ is the unique primitive matrix in its homology class.
\qed\\\\

\section{Proof of Theorem \ref{primitivesignedtheorem}}

 We will follow Henrich's proof of Theorem \ref{Allison}.  Let $(P, s,d, b,\epsilon)$ and $(P',s',d',b',\epsilon')$ be two homologous primitive signed singular based matrices.  We will show that the sequence of $M''_i$, $(M'')^{-1}_i$, and $N''$ moves relating two primitive singular based matrices can be replaced by a sequence of moves of the form $A\circ B \circ I \circ C$, where $A$ is a composition of extensions and $N''$ moves, $B$ is a single $D''_{ij}$ or $N''$ move, $I$ is an isomorphism, and $C$ is a composition of inverse extensions and $N''$ moves.  Since one cannot apply inverse extensions to a primitive matrix or obtain a primitive matrix by applying extensions to another matrix, the sequence must be of the form $B\circ I$, and the theorem follows.\\\\
As in Henrich's proof, we will show:
\begin{itemize}
\item The $D''_{ij}$ moves commute with inverse extensions and extensions,
\item A sequence of $D''_{ij}$ and $N''$ moves can be replaced with a sequence containing at most one $D''_{ij}$ or $N''$ move (and possibly an isomorphism), and
\item A sequence of extensions, inverse extensions, and $N''$ moves can be rewritten so that all inverse extensions occur before all extensions.
\end{itemize}
These three claims allow us to put our sequence in the form $A\circ B \circ I \circ C$ above.  To see why, consider the leftmost $D''_{ij}$ move (if such a move exists).  By the first bullet above, we slide this move past extensions and inverse extensions until it is adjacent to a sequence of $N''$ and $D''_{ij}$ moves.  Now by the second bullet, replace this sequence of our leftmost $D''_{ij}$ move and other $D''_{ij}$ and $N''$ moves by a sequence containing a single $D''_{ij}$ or $N''$ move.  Now find the new leftmost $D''_{ij}$ move, and repeat this process.  In the end we will have a sequence of extensions, inverse extensions, and $N''$ moves possibly followed by a $D''_{ij}$ move.  By the third bullet, we can reorder this sequence so that it is of the form: $(M''_i \text{ and } N'' \text{ moves})\circ ((M'')^{-1}_i \text{ and }N'' \text{ moves})\circ(\text{possibly a single } D''_{ij} \text{ move})$.  If there is a $D''_{ij}$ move at the end, we can move it between the inverse and ordinary extensions by sliding it past inverse extensions, and whenever it becomes adjacent to an $N''$ move, replace the result by a sequence containing a single $N''$ move or $D''_{ij}$ move as necessary.  

Our first step is to show that the $D''_{ij}$ moves commute with inverse extensions and extensions.  For $D''_{12}$ and $D''_{21}$, the proof is the same as Henrich's, so we only consider $D''_{33}$.  But this is clear because $D''_{33}$ just changes the sign of the matrix when the distinguished element is self-complementary, and a sign switch certainly commutes with inverse extensions and extensions.\\\\
Next we show that any sequence of $D''_{ij}$ and $N''$ moves can be replaced by a single $D''_{ij}$ or $N''$ move, or an isomorphism.  Henrich shows that this is true for compositions of $D''_{12}$, $D''_{21}$ and $N''$.  So we consider the compositions $D''_{33}\circ N''$, $N''\circ D''_{33}$, $D''_{33}\circ D''_{ij}$, and $D''_{ij}\circ D''_{33}$.  If the composition $D''_{33}\circ N''$ occurs, then the distinguished element must be self-complementary, implying that $N''$ is equivalent to a sign switch composed with an isomorphism.  Therefore the composition $D''_{33}\circ N''$ is an isomorphism.  This holds for $N''\circ D''_{33}$ as well.  The composition $D''_{33}\circ D''_{33}$ is also an isomorphism.  The compositions $D''_{33}\circ D''_{12}$, $D''_{33}\circ D''_{21}$, $D''_{12}\circ D''_{33}$, and $D''_{21}\circ D''_{33}$ can only occur if the core or annihilating-like distinguished element is also self-complementary.   Suppose that a core element $d$ is self-complementary.  Then $b(d,h)+b(d,h)=b(s,h)$ for all $h\in G$, but $b(d,h)=b(s,h)$, so $b(d,h)=0$.  Therefore $s$ is annihilating-like.  Similarly, if an annihilating element is self-complementry, $s$ is annihilating like.  If $s$ is annihilating like, the moves $D''_{12}$, $D''_{21}$, and $D''_{33}$ are all sign switches, so the composition of any two of them is an isomorphism.\\\\
Finally we must show that a sequence containing extensions, inverse extensions, and $N''$ moves, but no $D''_{ij}$ moves, can be rewritten so that inverse extensions occur before extensions.  Turaev \cite{Turaev2} showed that sequences containing only extensions and inverse extensions can be rewritten so that inverse extensions occur before extensions.  Like Henrich, we must consider sequences of the form $(M''_j)^{-1}\circ N'' \circ M''_i$ which are not equivalent to $D''_{ij}$ moves or isomorphisms.  In her case, such a sequence could be equivalent to an isomorphism, but in our case, this cannot happen because all of these sequences change the sign of the matrix exactly once.  In each case, we only need to consider the case where the $N''$ move interacts with both the extension and inverse extension, because if it does not, then the $N''$ move commutes with at least one of those moves, and then Turaev's results imply the sequence can be rewritten in the desired form.\\\\
{\bf Case 1: i=j=1.}  Henrich shows that this is equivalent to $D''_{12}=D''_{21}$, where $s$ is annihilating-like.\\
{\bf Case 2: i=j=2.} This is also equivalent to $D''_{12}=D''_{21}$, where $s$ is annihilating-like.\\ 
{\bf Case 3: i=1, j=2.}  This is $D''_{21}$.\\
{\bf Case 4: i=2, j=1.} This is $D''_{12}$.\\
{\bf Case 5: i=1, j=3.}  Henrich shows this is equivalent to $(M''_2)^{-1}\circ N''$.\\
{\bf Case 6: i=2, j=3.}  This is equivalent to $(M''_1)^{-1}\circ N''$.\\
{\bf Case 7: i=3, j=1.}  Henrich shows this is equivalent to $M''_1\circ D''_{12}$.\\
{\bf Case 8: i=3, j=2.}  This is equivalent to $M''_2\circ D''_{21}$.\\
{\bf Case 9: i=3, j=3.}  This is the case which differs from that of Henrich.  First $M''_3$ adds a pair of complementary elements $c_1$ and $c_2$.  The distinguished element $d$ must be complementary to one of these in order to apply $N''$.  So suppose $d$ and $c_1$ are complementary, so that $d$ and $c_2$ are the same with respect to $b$.  Then $N''$ makes $c_1$ the new distinguished element.  At the final stage, when we apply $(M''_3)^{-1}$, we have two options.  The first is that $d$ (and also both $c_i$) is self-complementary, so that $(M''_3)^{-1}$ removes $d$ and $c_2$.  In this case, the composition is a $D''_{33}$ move.  (Henrich did not need to consider this case because it is an isomorphism of unsigned matrices).  The other possibility is that $(M''_3)^{-1}$ only removes one of $d$ or $c_2$.  Since $b$ agrees on these elements, it does not matter which is removed.  Henrich shows that in this case, the composition is equivalent to an $N''$ move.
\qed

\section{Cases when the operation $\mu$ gives a formula for the minimal self-intersection number, and the proof that $\mu$ gives a stronger bound on $m([\alpha])$ than Turaev's cobracket $\nu$}\label{sec:bounds}
In this section, we prove Theorem \ref{primitiveformulaintro} and Corollary \ref{primitiveformulaintrocor}.  Recall that Theorem \ref{primitiveformulaintro} describes cases when the operation $\mu$ gives a formula for the minimal self-intersection number.  Corollary \ref{primitiveformulaintrocor} states that the bound on the minimal self-intersection number given by $\mu$ is stronger than the bound given by Turaev's cobracket $\nu$.

We are interested in describing strings $\alpha$ such that $m([\alpha])=t(\mu([\alpha]))/2 +n-1$, where $n\geq 1$ is the largest integer such that a minimal representative $\beta$ of $[\alpha]$ can be realized as a curve $B$ on an oriented surface $F$, and $\langle B\rangle=\langle \gamma \rangle ^n \in \pi_1(F)$.  Recall that by Proposition \ref{powers}, we always have $m([\alpha])\geq t(\mu([\alpha]))/2 +n-1$.

Theorem \ref{primitiveformulaintro} gives examples of strings $\alpha$ such that $m([\alpha])=t(\mu(\alpha))/2$.  Hence for these strings, $n=1$.

\begin{prop} \label{semitrivial} Suppose the term $[\alpha^\epsilon_e]$ of $\mu$ is semi-trivial.  Then the distinguished element of any primitive matrix in the homology class of $T(\alpha^\epsilon_e)$ is core or annihilating.
\end{prop}
\pp
If $[\alpha^\epsilon_e]$ is semi-trivial, then there exists some signed singular string $\tau^\delta_d$ in its class such that one of the arcs bounded by the endpoints $a$ and $b$ of $d$ contains no endpoints on its interior.  If we form the matrix $T(\tau^\delta_d)$, the distinguished element will be core or annihilating.  When the matrix $T(\tau^\delta_d)$ is reduced to a primitive matrix, the distinguished element can only change during an $N''$ move, so the distinguished element of the primitive matrix will be core or annihilating.
\qed

\ref{primitiveformulaintro}{ \bf Theorem. }{\it  Let $\alpha$ be a virtual string, whose based matrix $T(\alpha)$ is primitive and does not contain a self-complementary element.  Then $m([\alpha])=t(\mu([\alpha]))/2$.    If $T(\alpha)$ is primitive and does contain a self-complementary element, then either $m([\alpha])=t(\mu([\alpha]))/2$ or $m([\alpha])=t(\mu([\alpha]))/2+1$ .}

\pp  By Proposition \ref{semitrivial}, and the fact that $T(\alpha)$ is primitive, no term of $\mu$ is semitrivial.  Hence every term is nonzero.  Since $T(\alpha)$ is primitive (i.e., $T(\alpha)=T_\bullet(\alpha)$), every $e$ in $T(\alpha)$ is an element of $P(R)$, where $R$ is an empty sequence of moves (see Corollary \ref{uniquematrix} for the definition of $P(R)$).  

First we suppose $e$ is not self-complementary.  We will show that the term $[\alpha^+_e]$ cannot cancel with any other term $-[\alpha^-_f]$ of $\mu$ (and then , by symmetry, $-[\alpha^-_e]$ cannot cancel with $[\alpha^+_f]$).  Suppose that these terms cancel.  Then the signed singular based matrices $T(\alpha^+_e)$ and $T(\alpha^-_f)$ are homologous by Proposition \ref{homotopicimplieshomologoussignedsingular}.  By Corollary \ref{uniquematrix}, there is a unique primitive matrix $T_\bullet(\alpha^+_e)$ in the homology class of $T(\alpha^+_e)$, which is just the signed singular matrix with $T(\alpha)$ as the underlying based matrix, distinguished element $e$ and sign $+$.  Similarly, if $f$ is not self-complementary, we can find the unique primitive matrix $T_\bullet(\alpha^-_f)$ in the homology class of $T(\alpha^-_f)$ by using $T(\alpha)$ as the underlying based matrix, with distinguished element $f$ and sign $-$.  The matrices $T_\bullet(\alpha^+_e)$ and $T_\bullet(\alpha^-_f)$ have opposite signs, so they are not isomorphic.  If $f$ is not self-complementary, then we have reached a contradiction, so these terms cannot cancel.  If $f$ is self-complementary, the terms still cannot cancel.  By Theorem \ref{primitivesignedtheorem}, and because $e$ is not self-complementary and $f$ is, the matrices $T_\bullet(\alpha^-_f)$ and $T_\bullet(\alpha^+_e)$ are not related by the moves $N''$ or $D''_{ij}$.  

So the only way that two terms $[\alpha^+_e]$ and $-[\alpha^-_f]$ can cancel is if both $e$ and $f$ are self-complementary.  There is at most one self-complementary element in $T(\alpha)$, so $e=f$, and $[\alpha^+_e]$ might cancel with $-[\alpha^-_e]$, reducing the total number of terms by 2.
\qed\\\\

\subsection{Examples of classes $\alpha$ such that $\mu$ gives a formula for $m([\alpha])$.}  The purpose of Theorem \ref{primitiveformulaintro} is to find examples of classes such that $m([\alpha])=t(\mu(\alpha))/2+n-1$, in the case where $n=1$.  We now show that such classes exist.

Consider the string $\alpha_{p,q}$ defined in \cite[p. 2464]{Turaev2}.  This string is a copy of $S^1$ oriented counterclockwise, with $p$ vertical arrows pointing upward and $q$ horizontal arrows pointing from right to left, so that the copy of $S^1$ can be partitioned into four disjoint arcs containing the heads of the vertical arrows, the heads of the horizontal arrows, the tails of the vertical arrows, and the tails of the horizontal arrows, respectively.

Turaev finds a formula for the based matrix $T(\alpha_{p,q})$ (part of the formula is given in the next paragraph), and shows that if $p,q \geq 1$ and one of $p$ or $q$ is at least 2, then this matrix is primitive \cite[pp. 2464, 2471]{Turaev2}.  We will use this formula to check that $T([\alpha_{p,q}])$ does not contain a self-complementary element.  It will follow that for $p,q\geq 1$ and one of $p$ or $q \geq 2$, we have $m([\alpha_{p,q}])=t(\mu([\alpha_{p,q}]))/2=p+q$. 

Now we check that $T(\alpha_{p,q})$ does not contain a self-complementary element.  Label the vertical arrows $e_1, \dots , e_p$ from left to right, and label the horizontal arrows $e_{p+1}, \dots e_{p+q}$ from bottom to top. Then $b(e_i,s)=q$ for $i=1,\dots, p$ and $b(e_{p+j},s)=-p$ for $j=1,\dots, q$.  Also, $b$ vanishes on any pair of vertical arrows and $b$ vanishes on any pair of horizontal arrows.  We do not need to compute the rest of the matrix.  Suppose $c$ is a self-complementary element of $T(\alpha_{p,q})$.  Then $2b(g,c)=b(g,s)$ for all arrows $g$ of $\alpha_{p,q}$.  Suppose $c$ is horizontal.  Then choose some other horizontal arrow $h\neq c$.  We have $2b(h,c)=0$, but $b(h,s)=-p\neq 0$.  Hence $c$ cannot be horizontal.  Similarly it is easy to check $c$ cannot be vertical.  Thus $T(\alpha_{p,q})$ does not contain a self-complementary element, so by Theorem \ref{primitiveformulaintro}, we have $m([\alpha_{p,q}])=t(\mu([\alpha_{p,q}]))/2=p+q$.

 \ref{primitiveformulaintrocor} { \bf Corollary. }{\it The bound on $m([\alpha])$ given by $\mu$ is stronger than the bound given by Turaev's cobracket $\nu$.  Namely, the number of terms of $\mu([\alpha])$ is greater than or equal to the number of terms of $\nu[\alpha]$, and there are virtual homotopy classes $[\alpha]$ such that this inequality is strict.}

\pp  By Subsection \ref{sec.smoothingmap}, we have $\nu=S\circ \mu$.  It follows that $t(\mu([\alpha]))\geq t(\nu([\alpha]))$.   It is easy to check that $\nu([\alpha_{p,q}])=0$ for any $p,q$. On the other hand, we just saw that $m([\alpha_{p,q}])=t(\mu([\alpha_{p,q}]))/2=p+q$.  Hence $t(\mu([\alpha_{p,q}]))> t(\nu([\alpha_{p,q}]))$.

\qed

\section{An example showing the bound on $m([\alpha])$ given by $\mu$ is sometimes stronger than the bound $\rho$.}

In this section we prove that the string $M$ in Figure 1 satisfies $m([M])=|\arr(M)|=5$.  Then we show that $m([M])=t(\mu([M]))/2=5,$ while $\rho([M])=4.$

\subsection{Irreducible and minimal strings.}  We call a string $\alpha$ {\it crossing-reducible} if there is a string $\alpha'$ such that $\alpha$ and $\alpha'$ are related by a (possibly empty) sequence of Type 3 moves, and if a crossing-reducing Type 1 or Type 2 move can be applied to $\alpha'$.  A string which is not crossing-reducible is called {\it crossing-irreducible}.  A string $\alpha$ is {\it crossing-minimal} if $m([\alpha])=|\arr(\alpha)|$.

The notation $(F,\alpha)$ denotes an orientable surface $F$ with a curve $\alpha$ on it.  Recall that any virtual string $\alpha$ can be realized as a curve, which we also call $\alpha$, on a surface $F$.  There is a canonical realization of $\alpha$ on a surface, which can be obtained by gluing disks to the boundary components of the surface described in \cite[p. 2468]{Turaev2} (see also \cite{Carter}). Roughly speaking, one collapses each arrow of $\alpha$ to a point to get a framed $4$-valent graph, glues disks to each vertex, and then attaches bands between those disks along each edge in such a way that the resulting surface is orientable.   This canonical surface is the surface of smallest genus realizing $\alpha$.  We say two pairs $(F,\alpha)$ and $(F',\alpha')$ are virtually homotopic if the virtual strings realized by $\alpha$ and $\alpha'$ on $F$ and $F'$ respectively are virtually homotopic.

We call a pair $(F,\alpha)$ {\it genus-reducible} if there is a string $\alpha'$ homotopic to $\alpha$ on $F$, and a nontrivial simple closed curve $\gamma \subset F-\text{Im}(\alpha')$. If one cuts along $\gamma$ and glues disks to the resulting boundary components, the genus of the resulting surface is less than the genus of $F$; this process is called {\it destabilization} of $F$ along $\gamma$.  A pair which is not genus-reducible is called {\it genus-irreducible}.  We call a pair $(F,\alpha)$ {\it genus-minimal} if $F$ is the surface of smallest genus on which a representative of the virtual class $[\alpha]$ can be realized.

\subsection{The string $M$ is crossing-irreducible.}  Gibson showed, with the aid of a computer, that the string $M$ in Figure \ref{Gibsonex.fig} is crossing-irreducible \cite[p. 17, Table 6]{Gibson}.  (Gibson refers to this string by its nanoword $ABCADBECDE:bbbbb$).  

\subsection{The string $M$ is crossing-minimal.}  Kadokami stated that two crossing-irreducible strings are related by a (possibly empty) sequence of Type 3 moves \cite[Theorem 3.8]{Kadokami}.  It would follow from that statement that crossing-irreducible strings are crossing-minimal.  Kadokami proved this statement for flat virtual links, but Gibson unfortunately found a counterexample in the case where there is more than one component \cite[p. 18]{Gibson}.  However, the statement is true for virtual strings; it follows from the work of Ilyutko, Manturov and Nikonov \cite{IlyutkoManturovNikonov} as well as the work of Hass and Scott \cite{HassScott}, and we will show this below.  Then since Gibson showed $M$ is crossing-irreducible, Kadokami's statement will imply that the string $M$ is crossing-minimal.

We now explain why Kadokami's statement follows from \cite{IlyutkoManturovNikonov} and \cite{HassScott}. 

The results of Ilyutko, Manturov and Nikonov that we use are Theorem 3.2 and Corollary 3.1 in \cite{IlyutkoManturovNikonov}.  We state Theorem 3.2 in two ways.  The first statement (Theorem \ref{flatkup} below) is the same as the statement given in their work, but with our terminology and notation.

\begin{thm}[Ilyutko, Manturov, Nikonov] \label{flatkup} Let $(F,\alpha)$ and $(F',\alpha')$ be two virtually homotopic genus-minimal pairs.  Then there is a homeomorphism $\phi:F\rightarrow F'$ such that $\phi(\alpha)$ is homotopic to $\alpha'$ on $F'$.
\end{thm}
They actually prove the following stronger statement, Theorem \ref{flatkupstrong}, during their proof of Theorem 3.2 in \cite{IlyutkoManturovNikonov}.
\begin{thm}[Ilyutko, Manturov, Nikonov] \label{flatkupstrong} Let $(F,\alpha)$ and $(F',\alpha')$ be two virtually homotopic genus-irreducible pairs.  Then there is a homeomorphism $\phi:F\rightarrow F'$ such that $\phi(\alpha)$ is homotopic to $\alpha'$ on $F'$.
\end{thm}
Ilyutko, Manturov and Nikonov then deduce the following corollary using the results of Hass and Scott. We do not use this corollary, but state it here because it is similar to Kadokami's statement. 
\begin{cor}[Ilyutko, Manturov, Nikonov]\label{flatkupcor}  Let $\alpha$ and $\alpha'$ be two virtually homotopic crossing-minimal strings.  Then there is a (possibly empty) sequence of Type 3 moves taking $\alpha$ to $\alpha'$.    
\end{cor}

\begin{cor}[Kadokami's Statement]\label{flatkupcorstrong}  Let $\alpha$ and $\alpha'$ be two virtually homotopic crossing-irreducible strings.  Then there is a (possibly empty) sequence of Type 3 moves taking $\alpha$ to $\alpha'$.
\end{cor}

{\it Proof of Corollary \ref{flatkupcorstrong}}.  Suppose $\alpha$ and $\alpha'$ are crossing-irreducible.  Realize $\alpha$ and $\alpha'$ on their canonical surfaces $F$ and $F'$.  Hass and Scott \cite{HassScott} showed that any curve on a surface can be homotoped to a curve with minimal self-intersection without increasing the number of self-intersection points of the curve at any time during the homotopy, and in particular, two homotopic curves with minimal self-intersection are realted by Type 3 moves on the surface.  This implies that $\alpha$ and $\alpha'$ have minimal self-intersection in their free homotopy classes on $F$ and $F'$ respectively.  Using the result of Hass and Scott, along with the fact that $F$ and $F'$ are the canonical surfaces for $\alpha$ and $\alpha'$, it is straightforward to check that the pairs $(F,\alpha)$ and $(F', \alpha')$ are genus-irreducible; this argument is similar to an argument in the proof of Theorem \ref{flatkup}.  Hence by Theorem \ref{flatkupstrong} there is a homeomorphism $\phi:F\rightarrow F'$ such that $\phi(\alpha)$ is homotopic to $\alpha'$ on $F'$.    Since $\phi(\alpha)$ and $\alpha'$ are crossing-irreducible as virtual strings, they must both have the fewest number of self-intersection points of any curve in their free homotopy class on $F'$.  Therefore either $\phi(\alpha)$ and $\alpha'$ are related by a regular isotopy, or there is a sequence of Type 3 moves taking $\phi(\alpha)$ to $\alpha'$ on $F'$.  Hence $\alpha$ and $\alpha'$ are related by a (possibly empty) sequence of Type 3 moves. \qed

Now that we have established Corollary \ref{flatkupcorstrong}, we know that $M$ is crossing-minimal.

 \subsection{The proof that $m([M])=5=t(\mu([M]))/2$.}  In Lemma \ref{minimalnotsemitrivial}, we will show that since $M$ is crossing-minimal, no term of $\mu([M])$ is semi-trivial.  Now in order to show that $t(\mu([M]))/2=5$, we just need to show that no two terms of $\mu([M])$ cancel with each other.  The signed singular matrices $T(M^+_A)$, $T(M^+_B)$, $T(M^+_D)$, $T(M^+_E)$, $T(M^-_A)$, $T(M^-_B)$, $T(M^-_D)$, and $T(M^-_E)$ become primitive after the annihlating element $C$ is removed.  The resulting primitive matrices are the unique primitive matrices in their homology classes, and they are all distinct.  The matrices $T(M^+_C)$ and $T(M^-_C)$ are primitive, and are different than those mentioned above because they contain 6 elements rather than 5.  They are also not homologous to each other because they do not differ by a $D''_{12}$, $D''_{21}$, $D''_{33}$ or $N''$ move.  Hence the five signed singular homology classes corresponding to the five positive terms of $\mu([M])$ are different from the five signed singular homology classes corresponding to the five negative terms of $\mu([M])$, so no two terms of $\mu([M])$ cancel.

\begin{figure}[htbp]
   \centering
   \includegraphics[width=3.5cm]{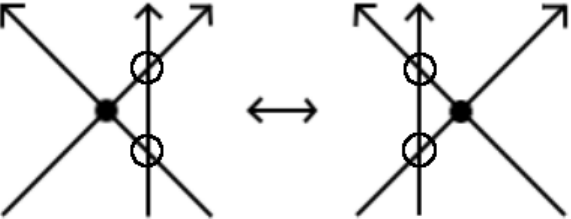} 
  \caption{A move for planar diagrams of singular virtual strings with one of several possible choices of orientation on the branches.}
   \label{type3movewithdotvirtual.fig}
\end{figure}

In the proof of Lemma \ref{minimalnotsemitrivial}, it will be useful to view singular virtual strings as flat singular virtual knot diagrams in the plane, up to the virtual diagram moves in \cite{Henrich}.  These moves are flat versions of the usual virtual Reidemeister moves in \cite{Kauffman}, plus the moves in Figures \ref{Type2movewithdot.fig} and \ref{Type3movewithdot.fig} (ignoring signs) and the move in Figure \ref{type3movewithdotvirtual.fig}.
\begin{lem}  \label{minimalnotsemitrivial}Suppose $\alpha$ is a virtual string such that $| \arr(\alpha)|=m([\alpha])$.  Then none of the singular strings $\alpha_e$, $e\in \arr(\alpha)$, are semi-trivial.
\end{lem}
\pp  Let $\alpha$ be a virtual string such that $|\arr(\alpha)|=m([\alpha])$, and suppose $\alpha_e$ is semi-trivial for some $e\in \arr(\alpha)$.  Then $\alpha_e$ is homotopic to a singular string $\tau_f$ where $f=(a,b)$ and either the arc $ab$ or the arc $ba$ contains no endpoints of arrows of $\tau$.

\begin{figure}[htbp]
   \centering
   \includegraphics[width=3.5cm]{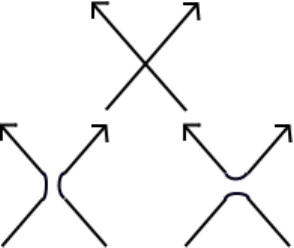} 
  \caption{Smoothing a crossing of $\alpha$ with and against the orientation of $\alpha$.}
   \label{smoothings.fig}
\end{figure}

We will see that we can smooth $\alpha$ at the crossing corresponding to the arrow $e$ in such a way that the resulting virtual string $\alpha'$ is in the same virtual homotopy class as $\alpha$.  This string has one less crossing than $\alpha$, contradicting the assumption that $|\arr(\alpha)|=m([\alpha])$.

\begin{figure}[htbp]
   \centering
   \includegraphics[width=3.5cm]{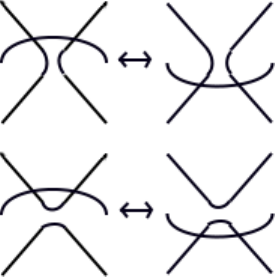} 
  \caption{Type 3 moves with center crossing smoothed.}
   \label{type3smoothed.fig}
\end{figure}

\begin{figure}[htbp]
   \centering
   \includegraphics[width=3.5cm]{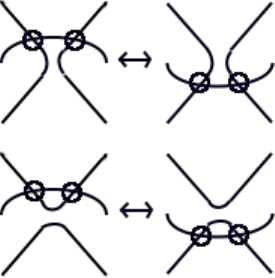} 
  \caption{Virtual Type 3 moves with center crossing smoothed.}
   \label{type3smoothedvirtual.fig}
\end{figure}

\begin{figure}[htbp]
   \centering
   \includegraphics[width=3cm]{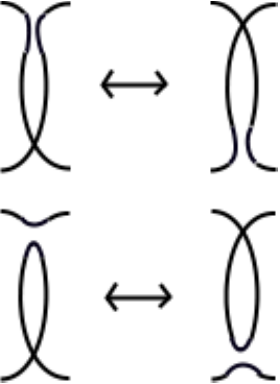} 
  \caption{Virtual Type 2 move with center crossing smoothed.}
   \label{type2smoothed.fig}
\end{figure}

Now we construct $\alpha'$ and show why it is homotopic to $\alpha$.  First smooth $\alpha$ at its distinguished crossing $e$ against the orientation of $\alpha$ (see Figure \ref{smoothings.fig}).  This is the virtual string $\alpha'$, though we have not yet specified its orientation.  The string $\alpha'$ is homotopic to the string $\tau$, where $\tau$ is the ordinary string underlying the singular string $\tau_f$.  The homotopy from $\alpha'$ to $\tau$ is given by using the sequence of moves for singular strings taking $\alpha_e$ to $\tau_f$, but where the moves involving the distinguished crossing in Figures \ref{Type3movewithdot.fig}, \ref{Type2movewithdot.fig}, and \ref{type3movewithdotvirtual.fig} are replaced with moves in Figures \ref{type3smoothed.fig}, \ref{type2smoothed.fig}, and \ref{type3smoothedvirtual.fig} respectively, followed by a single Type 1 move, which adds the arrow $f$.  We pick an orientation of $\alpha'$ so that this resulting virtual string is the string $\tau$.  Now $\alpha'$ is homotopic to $\tau$, and $\tau$ is homotopic to $\alpha$, so $\alpha'$ must be homotopic to $\alpha$.  
\qed

\section{Comparing the bound given by $\mu$ to the bound $\rho$}  \label{sec.boundcomparison}
In this section we prove Theorem \ref{main2}.  Unlike Theorem \ref{primitiveformulaintro}, which only considers virtual classes that contain a representative $\alpha$ whose based matrix is primitive, Theorem \ref{main2} describes the bound on $m([\alpha])$ given by $\mu$ for an arbitrary class $[\alpha]$, and compares this bound to $\rho$.

\ref{main2} {\bf Theorem. } {\it Let $\alpha$ be any virtual string.  Then 
	$$m([\alpha])\geq t(\mu([\alpha]))/2\geq \rho([\alpha])-1+O/2.$$
	If the primitive based matrix associated to $[\alpha]$ does not contain a self-complementary element, then 
	$$m([\alpha])\geq t(\mu([\alpha]))/2\geq \rho([\alpha])+O/2.$$	
}

\subsection{Definitions of $O$, ordinary elements, positive and negative terms, and the standard primitive matrix}.   We call a term $\epsilon[\alpha^\epsilon_e]$ in the sum $\mu([\alpha])$ {\it positive} if its coefficient $\epsilon$ is positive and {\it negative} if its coefficient $\epsilon$ is negative.  If a term is semi-trivial, we put $\epsilon=0$ and do not consider the term to be positive or negative.  Note that the coefficient can be different from the sign of a primitive matrix in the homology class of $T(\alpha^\epsilon_e)$.
 
We call a core or annihilating element $e\in G-\{s\}$ of a based matrix $T=(G,s,b)$ {\it ordinary} if $s$ is not annihilating-like in $T$ (or equivalently, if the element $e$ is not self complementary, by Proposition \ref{annihilatinglike}).  If a signed singular homology class contains a primitive matrix with an ordinary core or annihilating distinguished element, then there are exactly two primitive matrices in its homology class which differ by a $D''_{12}$ or $D''_{21}$ move, by Theorem \ref{primitivesignedtheorem}.  Since these moves change the sign of the matrix, there is a {\it unique} primitive matrix in the homology class of $T(\alpha^\epsilon_e)$ with a {\it given sign} $\epsilon$.  We call this choice of primitive matrix the {\it standard} choice.  

Let $C$ (respectively, $A$) be the total number of positive terms $+[\alpha^+_e]$ such that the distinguished element in the standard primitive matrix in the homology class of $T(\alpha^+_e)$ is an ordinary core (respectively, annihilating) element.   Now put $O=|C-A|$.

\subsection{Proof of Theorem \ref{main2}.} First we prove the following lemma.  The notation $P(R)$ was introduced in Subsection \ref{sec.unique}.

\begin{lem}  \label{extra} Let $e\in P(R)$ for some fixed reduction $R$ of $T(\alpha)$, and assume $e$ is not self-complementary.  Suppose that there exists $f\notin P(R)$ so that the (unique) primitive signed singular based matrix corresponding to $[\alpha^+_e]$ is isomorphic to the (unique) primitive signed singular based matrix corresponding to $[\alpha^-_f]$.  Then there exists some $f'\notin P(R)$ such that $f'\neq f$ and the primitive signed singular based matrix corresponding to $[\alpha^+_{f'}]$ is isomorphic to those corresponding to $[\alpha^+_e]$ and $[\alpha^-_f]$.
\end{lem}
\pp  Suppose $e$ and $f$ satisfy the hypotheses in the statement of Lemma \ref{extra}.  Form the matrices $T(\alpha^+_e)$ and $T(\alpha^-_f)$.  These matrices have $T(\alpha)$ as their underlying nonsingular based matrix.  $R$ is a sequence of inverse extensions which can be applied to $T(\alpha)$, so the extensions in this sequence can be applied to $T(\alpha^+_e)$ and $T(\alpha^-_f)$ as long as they do not remove the distinguished element.  So begin reducing $T(\alpha^+_e)$ and $T(\alpha^-_f)$ according to the moves in $R$ until this is no longer possible, i.e., until just before a move in $R$ would remove the distinguished element.  Since $e\in P(R)$, the reduction $R$ will reduce $T(\alpha^+_e)$ completely, giving us a primitive matrix.  Since $f\notin P(R)$, the reduction $R$ of $T(\alpha^-_f)$ will at some point require the removal of $f$, so we stop the reduction.  Let $M_f$ be the signed singular based matrix at this stage.  At this point, the distinguished element $f$ is core, annihilating, or part of a complementary pair.  Notice that $f$ cannot be core or annihilating, because we are assuming that the primitive matrix obtained by reducing $T(\alpha^-_f)$ is isomorphic to that obtained by reducing $T(\alpha^+_e)$; if at some point during the reduction $f$ becomes core or annihilating, then the distinguished element in the primitive matrix must also be core or annihilating, and that is impossible.  So $f$ is part of at least one complementary pair.  The reduction $R$ is about to remove $f$ with some specific choice of element $f'$ complementary to $f$.  

Now consider the term $[\alpha^+_{f'}]$.  Form its matrix $T(\alpha^+_{f'})$.  Begin reducing $T(\alpha^+_{f'})$ according to the sequence $R$ of inverse extensions.  At some point, $R$ removes $f'$, but we already know this happens at the same point when $R$ removes $f$, i.e., they are removed together as part of a complementary pair.  So at this point, we can apply an $N''$ move to switch the distinguished element from $f'$ to $f$.  After we apply $N''$ we obtain the matrix $M_f$.  Thus we have partial reductions of the matrices $T(\alpha^-_f)$ and $T(\alpha^+_{f'})$ such that at some stage in the reduction, the matrices are identical.  Thus the primitive signed singular based matrices in the homology classes of $T(\alpha^-_f)$ and $T(\alpha^+_{f'})$ are isomorphic.  

\qed \\\\

{\it Proof of Theorem \ref{main2}.} We break the proof into two steps.

In the first step, we prove:

{\it  The number of terms of $\mu([\alpha])$ is at least $2\rho(\alpha)-2$.  If no element of $P(R)$ is self-complementary for some fixed reduction $R$ of $T(\alpha)$, then the number of terms of $\mu([\alpha])$ 
is at least $2\rho([\alpha])$.}

Proof of Step 1: We will show that the arrows $e\in P(R)$ which are not self-complementary in $T_\bullet(\alpha)$ each contribute two terms to $\mu([\alpha])$ that do not cancel with any other terms.  By Proposition \ref{semitrivial}, if $e\in P(R)$, then the terms $[\alpha^+_e]$ and $-[\alpha^-_e]$ are not semi-trivial.\\\\

 To prove the theorem, we will partition the terms of $\mu$ into sets such that two terms of $\mu$ are in the same set if and only if they have same primitive signed singular matrix, and count the number of positive and negative terms in each set. \\\\
We first consider all of the positive terms of $\mu$ of the form $[\alpha^+_e]$ where $e\in P(R)$ is not self-complementary.  Let $I^+_1,I^+_2,\dots, I^+_n$ be sets containing these terms, such that two terms are in the same set if and only if they have isomorphic primitive signed singular based matrices.  (Recall that the primitive matrix is unique when $e\in P(R)$ and is not self-complementary by Corollary \ref{uniquematrix}).  Similarly let $I^-_1, I^-_2,\dots, I^-_m$ be sets containing the negative terms $-[\alpha^-_e]$ where $e\in P(R)$ is not self-complementary, such that two terms are in the same set if and only if they have isomorphic primitive signed singular based matrices.  By Corollary \ref{uniquematrix}, the terms $[\alpha^+_e]$ and $-[\alpha^-_f]$ corresponding to non-self-complementary elements $e, f\in P(R)$ cannot cancel with other such terms.  This implies that the sets $I^+_k$ are disjoint from the sets $I^-_k$.  

Now we assign the terms $\epsilon[\alpha^\epsilon_f]$ for all $f\notin P(R)$ to the sets $I^\pm_k$ as follows:  Consider a primitive signed singular matrix $P_f$ homologous to $T(\alpha^\epsilon_f).$  If $P_f$ is isomorphic to the primitive matrix of the term $\delta[\alpha^\delta_e]$ in some $I^\delta_k$, then add $\epsilon[\alpha^\epsilon_f]$ to $I^\delta_k$.  By Lemma \ref{extra}, if $P_f$ is isomorphic to the primitive matrix of some term $\delta[\alpha^\delta_e]$, then there is some unique $f'\notin P(R)$ determined by $R$ such that $P_f$ is also isomorphic to the primitive matrix of a term $-\epsilon[\alpha^{-\epsilon}_{f'}]$. So both $\epsilon[\alpha^\epsilon_f]$ and $-\epsilon[\alpha^{-\epsilon}_{f'}]$ are in the set $I^\delta_k$, and one is a positive term while the other is a negative term.\\\\
Now each $I^\delta_k$ contains an equal number of positive and negative terms corresponding to arrows not in $P(R)$, as well as a certain number $n^\delta_k$ of terms corresponding to arrows in $P(R)$.  This means that there are at least $n^\delta_k$ terms in $I^\delta_k$ which do not cancel with any other terms of $\mu$.  Note that the sum of all the $n^\delta_k$ is either $2\rho(\alpha)$ if there are no self-complementary elements in $P(R)$ or $2\rho(\alpha)-2$ if there is a self-complementary element in $P(R)$, which concludes the proof of Step 1.

Step 2:  We prove the following statement :

{\it  The number of terms of $\mu([\alpha])$ is at least $2\rho(\alpha)-2+O$.  If the primitive matrix $T(\alpha)$ does not contain any self-complementary elements, then the number of terms is at least $2\rho(\alpha)+O$.}

{Proof of Step 2:}   The number of positive and negative terms with ordinary core (respectively, annihilating) distinguished elements in their standard primitive matrices are equal. A positive (respectively, negative) term with an ordinary core distinguished element in its standard primitive matrix  can only cancel with a negative (respectively, positive) term with an ordinary annihilating distinguished element in its standard primitive matrix.  Even if all terms satisfying the condition in the previous sentence cancel with each other, there will still be $|C-A|$ such terms leftover that do not cancel with any terms of $\mu$.
\qed\\\\

{\bf Acknowledgements:} This paper was revised while visiting the Max Planck Institute for Mathematics in Bonn, Germany, and I would like to thank the institute for its hospitality.  
I am very grateful to Vladimir Chernov for many helpful discussions, for providing translations of some of the references, and for commenting on drafts of this paper.  I would also like to thank V. Manturov for many helpful discussions concerning algorithms for finding minimal representatives of virtual homotopy classes, and to V. Manturov and D. Ilyutko for helpful discussions about their proof of the Kuperberg Theorem for flat virtual knots in \cite{IlyutkoManturovNikonov}.  Finally I would like to thank M. Chas for informing me about reference \cite{HassScott}.


\begin{thebibliography}{99999}

\bibitem{AMR2}
{\bf J.E. Andersen, J. Mattes, N. Reshetikhin}, {\it Quantization of the algebra of chord diagrams,} Math. Proc. Cambridge Philos. Soc., Vol. 124 no. 3 (1998), pp. 451-467

\bibitem{AMR}
{\bf J.E. Andersen, J. Mattes, N. Reshetikhin}, {\it The Poisson structure on the moduli space of flat connections and chord diagrams.} Topology 35 (1996), no. 4, 1069-1083.

\bibitem{Cahn}
{\bf P. Cahn}, {\it A Generalization of the Turaev Cobracket and the Minimal Self-Intersection Number}. New York J. Math 19 (2013), 253-283.

\bibitem{Carter}
{\bf S. Carter}, {\it Classifying immersed curves.}  Proc. Amer. Math. Soc. 111 (1991), no. 1, 281-287.

\bibitem{CarterKamadaSaito}
{\bf S. Carter, S. Kamada, and M. Saito}, {\it Stable equivalence of knots on surfaces and virtual knot cobordisms}.  J. Knot Theory Ramifications 11 (2002), 311-322.

\bibitem{Chas}
{\bf M. Chas}, {\it Combinatorial Lie Bialgebras of curves on surfaces}.  Topology 43 (2004), no. 3, 543-568.

\bibitem{Gibson}
{\bf A. Gibson}, {\it On tabulating virtual strings.} Acta Math. Vietnam. 33 (2008), no. 3, 493–518. arXiv:0808.0064v1.

\bibitem{HassScott}
{\bf J. Hass and P. Scott}, {\it Shortening Curves on Surfaces}.  Topology 33 (1994), no. 1, 25-43.

\bibitem{Henrich}
{\bf A. Henrich}, {\it A sequence of degree one Vassiliev invariants for virtual knots}.  Journal of Knot Theory and its Ramifications 19 (2010), no. 4, 461-487.


\bibitem{IlyutkoManturovNikonov}
{\bf D. P. Ilyutko, V.O. Manturov, I.M. Nikonov,} {\it Chetnost v teorii uzlov i graf-zaceplenij.}  To appear in Sovrem. Mat. Fundam. Napravl.

\bibitem{Kadokami}
{\bf T. Kadokami,} {\it Detecting Non-Triviality of Virtual Links.}  J. Knot Theory Ramifications 12 (2003), no. 6, 781-803.


\bibitem{Kauffman}
{\bf L. Kauffman,} {\it Virtual knot theory}.  Europ. J. Combinatorics 20 (1999), no. 7, 663-690.



\bibitem{Turaev2}
{\bf V. Turaev}, {\it Virtual strings,} Ann. Inst. Fourier (Grenoble) {\bf 54} (2004), no. 7, 2455-2525.

\bibitem{Turaev}
{\bf V. Turaev}, {\it Skein quantization of Poisson algebras of loops on surfaces}.  Ann. Sci. Ecole Norm. Sup.  (4) 24 (1991), no. 6, pp. $635-704$.

\bibitem{TuraevViro}
{\bf V. G. Turaev and O. Ya. Viro}, {\it Intersection of loops in two-dimensional manifolds.  II.  Free loops}.  Mat. Sbornik 121:3 (1983) $359-369$ (Russian); English translation in Soviet Math. Sbornik.
\end{thebibliography}
\end{document}